\documentclass[12pt]{article}
\usepackage{amsfonts}
\usepackage{amsfonts,mathrsfs}
\usepackage{CJK,fancyhdr,amscd}
\usepackage{anysize}
\usepackage{graphicx,epsfig}
\usepackage{amsthm,latexsym,amsmath,amssymb}
\usepackage[perpage,symbol*]{footmisc}

\newtheorem{thm}{Theorem}[section]

\newtheorem{prop}[thm]{Proposition}
\newtheorem{lem}[thm]{Lemma}

\newtheorem{rem}[thm]{Remark}
\def\pf{\noindent{\bf Proof.} }

\makeatletter \@addtoreset{equation}{section} \makeatother
\begin{document}
\baselineskip=20pt  \hoffset=-3cm \voffset=0cm \oddsidemargin=3.2cm
\evensidemargin=3.2cm \thispagestyle{empty}\vspace{10cm}
\parskip 8pt
 \vskip 2mm

\begin{center}
{\Large\bf Existence and multiplicity of Homoclinic solutions for
the second order Hamiltonian systems}
\end{center}

\centerline{\large Chungen
Liu $^{{\rm b},*,\ddag}$, Qingye Zhang $^{{\rm a},\dag}$ $\quad$ }
\begin{center}
$^{\rm a}${\it\small School of Mathematical Science, Nankai University\\
Tianjin 300071, P.R. China}
\end{center}

\begin{center}
$^{\rm b}${\it\small School of Mathematical Science and LPMC, Nankai
University\\Tianjin 300071, P.R. China}
\end{center}

\footnotetext[0]{$^*$Corresponding author}

\footnotetext[0]{\;{\it E-mail address}: zhzhy323@mail.nankai.edu.cn
(Qingye Zhang), liucg@nankai.edu.cn (Chungen Liu)}

\footnotetext[0]{$^\dag$Partially  supported by NFSC of
China(10701043).}

\footnotetext[0]{$^\ddag$Partially  supported by NFSC of
China(10531050,10621101) and 973 Program of STM(2006CB805903).}

\noindent{\bf Abstract}

{\small In this paper we study the existence and multiplicity of
homoclinic solutions for the second order Hamiltonian system
$\ddot{u}-L(t)u(t)+W_u(t,u)=0$, $\forall t\in\mathbb{R}$, by means
of the minmax arguments in the critical point theory, where $L(t)$
is unnecessary uniformly positively definite for all $t\in
\mathbb{R}$ and $W_u(t, u)$ sastisfies the asymptotically linear
condition.

\noindent{\it Keywords:} Homoclinic solution; Second order
Hamiltonian system; Linking structure; Mountain pass theorem

\noindent{\it MSC:} 37J45, 58E05, 34C37, 70H05}\\

\def\<{\langle}
\def\>{\rangle}

\section{Introduction and the main result}

Consider the second order Hamiltonian systems
\begin{equation}\label{HS}
\ddot{u}-L(t)u(t)+W_u(t,u)=0,\quad\forall t\in\mathbb{R},
\end{equation}
where $L\in C\left(\mathbb{R},\mathbb{R}^{N^2}\right)$ is a
symmetric matrix valued function, $W\in C^1\left(\mathbb{R}\times
\mathbb{R}^N , \mathbb{R}\right)$.   We say that a solution $u$ of
\eqref{HS} is homoclinic (to 0) if $u\in
C^2\left(\mathbb{R},\mathbb{R}^N\right)$, $u\neq 0$,
$u(t)\rightarrow 0$ and $\dot u(t)\rightarrow 0$ as $|t| \rightarrow
\infty$.

The existence and multiplicity of homoclinic solutions for
\eqref{HS} have been extensively investigated in many papers via the
variational methods, see, e.g., [1-6, 8, 9, 12, 14-19]. Most of them
treat the superquadratic case (see [1-6, 8, 9, 12, 14-16]), while
[18, 19] consider the asymptotically quadratic case and [5, 17]
treat the subquadratic case. But except for [5, 14] all known
results are obtained under the following assumption that $L(t)$ is
uniformly positively definite for all $t\in\mathbb{R}$, that is,
there exists a constant $l_0>0$ such that
\[
\<L(t)u, u\>\geq l_0|u|^2,\;t\in\mathbb{R},\;u\in\mathbb{R}^N,
\]
where $\<\cdot,\cdot\>$ and $|\cdot|$ are the  standard inner
product and the associated norm in $\mathbb{R}^N$ respectively and
we will always use these notations.

In this paper, we study the homoclinic solutions of \eqref{HS} where
$L(t)$ is unnecessary uniformly positively definite for all
$t\in\mathbb{R}$, and $W(t, u)$ satisfies subquadratic condition.
More precisely, $L$ satisfies

\noindent (L1) The smallest eigenvalue of $L(t)\rightarrow \infty$
as $|t|\rightarrow \infty$,  i.e.,
\[
l(t)\equiv \inf_{|u|=1,\ u\in \mathbb{R}^N }\<L(t)u, u\>\rightarrow
\infty,\quad\mbox{as $|t|\to\infty$},
\]

\noindent (L2) For some $a> 0$ and $\bar{r}> 0$, one of the
following is true:

$(i)$ $L\in C^1(\mathbb{R},\mathbb{R}^{N^2})$ and ${|L'(t)| \leq
a|L(t)|},\;\forall |t| \geq  \bar{r}$, or

$(ii)$ $L\in C^2(\mathbb{R},\mathbb{R}^{N^2})$ and $L''(t)\leq
aL(t),~\forall |t|\geq \bar{r}$,

\noindent where $L'(t)=(d/dt)L(t)$ and $L''(t)=(d^2/dt^2)L(t)$,

\noindent and $W(t, u)$ satisfies

\noindent (W1) $W(t,u)\geq 0$, $W(t,0)=0$ and $W_u(t,u)=o(|u|)$ as
$u\rightarrow 0$ uniformly in $t$, $|W_u(t,u)|\leq C_W(|u|)$ for
some $C_W>0$.

In what follows it will always be assumed that (L1) is satisfied.
Denote by A the selfadjoint extension of the operator
$-(d^2/dt^2)+L(t)$ with domain $D(A)\subset L^2\equiv
L^2\left(\mathbb{R},\mathbb{R}^N\right)$. Let
$\{E(\lambda):-\infty<\lambda<\infty\}$ and $|A|$ be the spectral
resolution and the absolute value of $A$ respectively, and
$|A|^{1/2}$ be the square root of $|A|$ with domain $D(|A|^{1/2})$.
Set $U=I-E(0)-E(-0)$, where $I$ is the identity map on $L^2$. Then
$U$ commutes with $A$, $|A|$ and $|A|^{1/2}$, and $A=U|A|$ is the
polar decomposition of $A$ (see \cite{K}). Let $E=D(|A|^{1/2})$, and
define on $E$ the inner product and norm
\begin{align*}
(u, v)_0&=\left(|A|^{1/2}u, |A|^{1/2}v\right)_{2}+(u,v)_{2},\\
\|u\|_0&=(u,u)_0^{1/2},
\end{align*}
where $(\cdot, \cdot)_{2}$ denotes the inner product in $L^2$; then
$E$ is a Hilbert space.

In order to learnt about the spectrum of $A$, We first need the
following lemma  from \cite{D} (cf. Lemma 2.1 in \cite{D}).

\begin{lem}\label{cptem}
Suppose that $L$ satisfies {\rm(L1)}, then $E$ is compactly embedded
in $L^p\equiv L^p(\mathbb{R}, \mathbb{R}^N)$ for $2\leq p\leq
\infty$.
\end{lem}

\begin{rem}\label{H12}
\rm It is easy to see that $E$ is continuously embedded in $H^{1,
2}(\mathbb{R}, \mathbb{R}^N)$ from the fact that
$C_0^\infty(\mathbb{R}, \mathbb{R}^N)$ is dense in $E$ and the proof
of Lemma 2.1 in \cite{D}.
\end{rem}

From \cite{D}, under the above assumption (L1) on $L$ and by Lemma
\ref{cptem}, we know that $A$ possesses a compact resolvent and the
spectrum $\sigma (A)$ consists of only eigenvalues numbered in
$\lambda_1\leq \lambda_2 \leq \cdots \to\infty,$ %
with a corresponding eigenfunctions $(e_n)(Ae_n=\lambda_ne_n)$,
forming an orthogonal basis in $L^2$. Let
$n^-=\#\{i|\lambda_i<0\}$, $n^0=\#\{i|\lambda_i=0\}$, %
and $\bar{n}=n^-+n^0$.
Set $E^-=\mathrm{span}\{e_1,\ldots, e_{n^-}\}$, %
$E^0=\mathrm{span}\{e_{n^-+1},\ldots, e_{\bar{n}}\}=\ker A$ %
and $E^+=\overline{\mathrm{span}\{e_{\bar{n}+1},\ldots\}}$. %
Then one has the orthogonal decomposition $E=E^-\oplus E^0\oplus
E^+$ with respect to the inner product $(\cdot, \cdot)_0$ on $E$.
Now we introduce on E the following inner product and norm:
\begin{align*}
   (u, v)&=\left(|A|^{1/2}u,|A|^{1/2}v\right)_{2}+\left(u^0,v^0\right)_{2},\\
  \|u\|&=(u,u)^{1/2},
\end{align*}
where $u=u^-+u^0+u^+$ and $v=v^-+v^0+v^+\in E=E^-\oplus E^0\oplus
E^+. $ Clearly the norms $\|\cdot\|$ and $\|\cdot\|_0$ are
equivalent (cf. \cite{D}). From now on $\|\cdot\|$ will be used.

\begin{rem}\label{orthdecom}
{\rm Note that the decomposition $E=E^-\oplus E^0\oplus E^+$ is also
orthogonal with respect to both $(\cdot, \cdot)$ and $(\cdot,
\cdot)_2$.}
\end{rem}

\begin{rem}
{\rm Since the norms $\|\cdot\|$ and $\|\cdot\|_0$ on $E$ are
equivalent, then by Lemma \ref{cptem}, for any $2\leq p\leq \infty$,
there exists $\beta_p>0$ such that
\begin{equation}\label{LpE}
|u|_p\leq \beta_p\|u\|,\;\forall u\in E,
\end{equation}
where $|\cdot|_p$ is the norm on $L^p$.}
\end{rem}

For later use, let
\begin{equation}\label{quadraticform}
a(u, v)=(|A|^{1/2}Uu, |A|^{1/2}v)_2,\quad \forall u, v\in E
\end{equation}
be the quadratic form associated with $A$. For any $u\in D(A)$ and
$v\in E$, we have
\begin{equation}\label{quadintform}
a(u, v)=\int\limits_\mathbb{R}\left(\<\dot{u}, \dot{v}\>+\<L(t)u,
v\>\right)dt
\end{equation}
and \eqref{quadintform} holds for all $u, v\in E$ since $D(A)$ is
dense in $E$. Moreover, by definition
\begin{equation}\label{quadnormform}
a(u, u)=((P^+-P^-)u,u)=\|u^+\|^2-\|u^-\|^2
\end{equation}
for all $u=u^-+u^0+u^+\in E$, where $P^\pm:E\to E^\pm$ are the
orthogonal projections with respect to the inner product $(\cdot,
\cdot)$.

We further make the following assumptions on $W$:

\noindent (W2) $W_u(t,u)=M(t)u +w_u(t,u)$ with $M$ a bounded,
continuous symmetric $N\times N$
 matrix-valued function and
$w_u(t,u)=o(|u|)$ as $|u|\rightarrow \infty$, $\forall
t\in\mathbb{R}$;

\noindent (W3) $m_0:= \displaystyle\inf_{t\in \mathbb{R}} \left[
\inf_{|u|=1,\ u\in \mathbb{R}^N }\<M(t)u, u\>\right]> \inf
\left(\sigma(A)\cap (0, \infty)\right)$;

\noindent (W4) $0\notin \sigma_p (A-M)$, where $\sigma_p (A-M)$ is
the point spectrum of $A-M$, $M$ is the operator defined on $L^2$ by
\[
(Mu)(t):=M(t)u(t), \;t\in \mathbb{R},\;u\in L^2.
\]

From the above spectral result of the operator A, the set
$\sigma(A)\cap(0, m_0)$ consists of only eigenvalues of finite
multiplicity, where $m_0$ is defined in (W3). Let $\ell$ denote the
number of  eigenvalues (counted with multiplicity) lying in $(0,
m_0)$.

Then we have our main result:

\begin{thm}\label{mainresult}
Suppose that {\rm (L1), (L2)} and {\rm (W1)\---(W4)} are satisfied.
Then \eqref{HS} has at least one nontrivial homoclinic solution. If
in addition $W(t,u)$ is even in $u$, then \eqref{HS} has at least
$\ell$ pairs of nontrivial homoclinic solutions.
\end{thm}
\newpage
\begin{rem}
\rm There are functions $L$ and $W$ which satisfy the conditions in
our Theorem \ref{mainresult} but do not satisfy the corresponding
conditions in [1-6, 8, 9, 12, 14-19]. For example, let
\[
  L(t) = \begin{cases}
(et^2-2)I_N,&   |t|\leq 1/\sqrt{e} \\
(\ln{t^2})I_N,&  |t|>1/\sqrt{e},
\end{cases}
\]
\[
W(t,
u)=\frac{1}{2}\left(e^{-t^2}+a\right)|u|^2\left(1-\frac{1}{\ln{(e+|u|)}}\right).
\]
Simple computation shows that $M(t)=\left(e^{-t^2}+a\right)I_N$ in
(W2) and we can choose suitable $a>\inf(\sigma(A)\cap (0, \infty))$
such that (W4) holds due to  the special spectral result of $A$
above.
\end{rem}

\section { Variational setting and proof of the main result}
 In order to establish a variational setting for the problem \eqref{HS},
we further need the following lemma which can be found in \cite{D}.

\begin{lem}[{\cite[Lemma 2.3]{D}}]\label{conem}
If L satisfies {\rm(L1)} and {\rm(L2)}, then $D(A)$ is continuously
embedded in $H^{2,2}\left(\mathbb{R},\mathbb{R}^N\right)$ and
consequently, we have
$$|u(t)|\rightarrow 0\;\mbox{\rm and} \;\dot u(t)\rightarrow 0\quad
\mbox{\rm as}\;|t| \rightarrow \infty,\;\;\forall u\in D(A).$$
\end{lem}
For any fixed $b>0$, let k be the number of eigenvalues of the
operator $A$(counted with multiplicity) lying in $[-b, b]$. Denote
by $f_i$ $(1\leq i\leq k)$ the corresponding eigenfunctions and set
$$L^{b-}:={\rm span}\{f_1,\ldots, f_k\},$$ then we have the orthogonal
decomposition $$L^2=L^{b-}\oplus L^{b+},\quad u=u^{b-} + u^{b+.}$$
where $L^{b+}$ is the orthogonal complement of $L^{b-}$ in $L^2.$
\\Correspondingly, $E$ has the decomposition $$E=E^{b-}\oplus
E^{b+}\;{\rm with}\; E^{b-}=L^{b-}\; {\rm and }\; E^{b+}=E\cap
L^{b+},$$ orthogonal with respect to both the inner products
$(\cdot, \cdot)_2$ and $(\cdot, \cdot)$. Then we have the following
lemma which will be used.

\begin{lem}\label{normcomp}
For any fixed $b>0,$ let $E=E^{b-}\oplus E^{b+}$ as above, then
\[
b|u|_2^2\leq \|u\|^2\quad{\rm for ~all}~u\in E^{b+},
\]
where $|\cdot|_2$ is the norm on $L^2$.
\end{lem}
\pf It is obvious from the definition of the norm $\|\cdot\|$ on $E$
and the distribution of the eigenvalues of $A$. $\hfill\Box$

By virtue of the quadratic form in \eqref{quadraticform}, we define
a functional $\it\Phi$ on $E$ by
\begin{align}\label{functional}
{\it\Phi}(u)&=\frac{1}{2}a(u,v)-{\it\Psi}(u)\notag\\[7pt]
&=\frac{1}{2}\int\limits_\mathbb{R}(|\dot{u}|^2+\<L(t)u,
u\>)dt-{\it\Psi}(u)\notag\\[7pt]
&=\frac{1}{2}\|u^+\|^2-\frac{1}{2}\|u^-\|^2-{\it\Psi}(u)\quad{\rm
where}\;{\it\Psi}(u)=\int\limits_\mathbb{R}W(t, u)dt
\end{align}
for all $u=u^-+u^0+u^+\in E=E^-\oplus E^0\oplus E^+ $. By (W1) and
Lemma \ref{cptem}, ${\it\Phi}$ and ${\it\Psi}$ are well defined.
Furthermore, we have
\begin{prop}\label{propofPhi}
Let $\rm(L1)$, $\rm(L2)$ and $\rm(W1)$ be satisfied. Then
${\it\Psi}\in C^1(E, \mathbb{R})$, and hence ${\it\Phi}\in C^1(E,
\mathbb{R})$. Moreover,
\begin{align}
{\it\Psi}'(u)v&=\int\limits_\mathbb{R}\< W_u(t, u),
v\>dt\label{derivativePsi}\\[5pt]
{\it\Phi}'(u)v&=(u^+, v^+)-(u^-, v^-)-{\it\Psi}'(u)v\notag\\[7pt]
&=(u^+, v^+)-(u^-, v^-)-\int\limits_\mathbb{R}\< W_u(t, u),
v\>dt\label{derivative}
\end{align}
for all $u=u^-+u^0+u^+\in E=E^-\oplus E^0\oplus E^+ $ and
$v=v^-+v^0+v^+\in E=E^-\oplus E^0\oplus E^+$, and critical points of
${\it\Phi}$ on $E$ are homoclinic solutions of \eqref{HS}.
\end{prop}
\pf We first verify \eqref{derivativePsi} by definition. Let $u\in
E$. Using (W1), by the mean value theorem and the H\"{o}lder
inequality, we have
\begin{align}\label{tgeqT}
&\left|\int\limits_{|t|>T}(W(t,u+v)-W(t,u)-\< W_u(t, u),
v\>)dt\right|\notag\\
&\leq C\left(\int\limits_{|t|>T}(|u|+|v|)^2dt\right)^{1/2}|v|_2\notag\\
&\leq
C\beta_2\left(\int\limits_{|t|>T}(|u|+|v|)^2dt\right)^{1/2}\|v\|,\;\forall\;T>0,\forall\;v\in
E,
\end{align}
where $C$ is a constant and the last inequality holds by
\eqref{LpE}. In view of Lemma \ref{cptem}, for any $\varepsilon>0$,
there is a $\delta_1>0$ and $T_\varepsilon>0$ such that
\begin{equation}\label{Cbeta2}
C\beta_2\left(\int\limits_{|t|>T_\varepsilon}(|u|+|v|)^2dt\right)^{1/2}\leq
\varepsilon/2,
\end{equation}
for all $v\in E$, $\|v\|\leq\delta_1$.\\
From Remark \ref{H12}, $u\in H^{1,2}(\mathbb{R}, \mathbb{R}^N)$.
Define ${\it\Psi}_T:E\to \mathbb{R}$ by
\[
{\it\Psi}_T(u):=\int\limits_{-T}^TW(t,u)dt,\;\forall u\in E.
\]
It is known (see, e.g.,\cite{R}) that ${\it\Psi}_T\in
C^1(H^{1,2}([-T, T], \mathbb{R}^N),\mathbb{R})$ for any  $T>0$.
Therefore, for the $\varepsilon$ and $T_\varepsilon$ given above, by
Remark \ref{H12}, there is a $\delta_2=\delta_2(\varepsilon,
T_\varepsilon, u)$ such that
\begin{equation}\label{tleqT}
\left|\int\limits_{-T_\varepsilon}^{T_\varepsilon}(W(t,u+v)-W(t,u)-\<
W_u(t, u),v\>)dt\right|%
\leq\frac{\varepsilon}{2}\|v\|,
\end{equation}
for all $v\in E$, $\|v\|\leq\delta_2$.\\
Combining  \eqref{tgeqT}, \eqref{Cbeta2} with \eqref{tleqT} and
taking $\delta=\min\{\delta_1,\delta_2\}$, then we obtain
\[
\left|\int\limits_\mathbb{R}(W(t,u+v)-W(t,u)-\< W_u(t, u),
v\>)dt\right|\leq \varepsilon \|v\|
\]
for all $v\in E$, $\|v\|\leq\delta$. Thus \eqref{derivativePsi}
follows immediately by the definition of Fr\'{e}chet derivatives.
Due to the form of ${\it\Phi}$ in \eqref{functional},
\eqref{derivative} also holds.

We then prove that ${\it\Psi}'$ is continuous. Suppose $u_n\to u_0$
in $E$ and hence $u_n\to u_0$ in $L^\infty$ by Lemma \ref{cptem}.
Note that
\begin{align}\label{Phidercont}
\sup\limits_{\|v\|=1}\left\|({\it\Psi}'(u_n)-{\it\Psi}'(u_0))v\right\|%
&=\sup\limits_{\|v\|=1}\left|\int\limits_\mathbb{R}\<W_u(t,u_n)-W_u(t,u_0),v\>dt\right|\notag\\
&\leq\sup\limits_{\|v\|=1}\left(\int\limits_\mathbb{R}|W_u(t,u_n)-W_u(t,u_0)|^2dt\right)^{1/2}|v|_2\notag\\
&\leq\beta_2\left(\int\limits_\mathbb{R}|W_u(t,u_n)-W_u(t,u_0)|^2dt\right)^{1/2}
\end{align}
where $\beta_2$ is the constant in \eqref{LpE}.\\
Note that, by Lemma \ref{cptem}, $(u_n)$ is bounded in $L^2$ since
$u_n\to u_0$ in $E$, i.e., there exists a constant $M_0>0$ such that
$|u_n|_2\leq M_0$, $\forall n\in \mathbb{N}$.  By (W1), for any
$\varepsilon>0$, there exists $\eta>0$ such that
\begin{equation}\label{Wueta}
|W_u(t, u)|\leq \frac{\varepsilon}{2(M_0+|u_0|_2)}|u|,\;\forall u\in
\mathbb{R}^N,\;|u|\leq\eta.
\end{equation}
Due to $u_0\in H^{1,2}(\mathbb{R}, \mathbb{R}^N)$ and $u_n\to u_0$
in $L^\infty$, there exist $T_\varepsilon>0$ and $N_1\in\mathbb{N}$
such that for all $n>N_1$ and $|t|\geq T_\varepsilon$, it holds that
\begin{equation}\label{ngeqN1}
\begin{aligned}
|W_u(t, u_n(t))|&\leq \frac{\varepsilon}{2(M_0+|u_0|_2)}|u_n(t)|,\\
|W_u(t, u_0(t))|&\leq \frac{\varepsilon}{2(M_0+|u_0|_2)}|u_0(t)|.
\end{aligned}
\end{equation}
Observe also that $(u_n)$ is bounded in $L^\infty$, then by (W1) and
Lebesgue's Dominated Convergence Theorem,
\[
\left(\int\limits_{-T_\epsilon}^{T_\varepsilon}|W_u(t,u_n)-W_u(t,u_0)|^2dt\right)^{1/2}\rightarrow
0,\;{\rm as}\;n\to \infty.
\]
Then there exists $N_2\in \mathbb{N}$ such that for all $n>N_2$,
\[
\left(\int\limits_{-T_\epsilon}^{T_\varepsilon}|W_u(t,u_n)-W_u(t,u_0)|^2dt\right)^{1/2}\leq
\varepsilon/2
\]
Combining this with \eqref{ngeqN1} and taking $N=\max\{N_1, N_2\}$,
we have
\begin{align*}
&\left(\int\limits_\mathbb{R}|W_u(t,u_n)-W_u(t,u_0)|^2dt\right)^{1/2}\\
&\leq
\left(\int\limits_{-T_\epsilon}^{T_\varepsilon}|W_u(t,u_n)-W_u(t,u_0)|^2dt\right)^{1/2}+%
\left(\int\limits_{|t|>T_\varepsilon}|W_u(t,u_n)-W_u(t,u_0)|^2dt\right)^{1/2}\\
&\leq\frac{\varepsilon}{2}+\frac{\varepsilon}{2(M_0+|u_0|_2)}(|u_n|_2+|u_0|_2)\leq\varepsilon
\end{align*}
for all $n>N$. This shows that
\[
\left(\int\limits_\mathbb{R}|W_u(t,u_n)-W_u(t,u_0)|^2dt\right)^{1/2}\to
0,\;n\to \infty.
\]
Thus the continuity of ${\it\Psi}'$ follows immediately by
\eqref{Phidercont}. Consequently, the form of ${\it\Phi}$ yields
${\it\Phi}\in C^1(E, \mathbb{R})$.

Finally, we show that critical points of ${\it\Phi}$ on $E$ are
homoclinic solutions of \eqref{HS}. Note first that, by means of a
standard argument, \eqref{quadraticform}---\eqref{quadnormform} and
\eqref{derivative} imply that critical points of ${\it\Phi}$ belong
to $C^2\left(\mathbb{R},\mathbb{R}^N\right)$ and satisfy \eqref{HS}.
Now for any critical point $u$ of ${\it\Phi}$ on $E$, by (W1) and
Lemma \ref{cptem}, one has
\begin{align*}
|Au|_2^2&=\int\limits_\mathbb{R}|W_u(t,u)|^2dt\\[5pt]
&\leq C_W^2|u|_2^2<\infty.
\end{align*}
where $C_W$ is the constant in (W1). Thus $u\in D(A)$ and $u$ is a
homoclinic solution of \eqref{HS} by Lemma \ref{conem}. The proof is
completed. $\hfill\Box$

We will make use of minimax arguments to prove our main result and
first state two results of this type from Rabinowitz \cite{R} and
Ghoussoub \cite{G} here. One is the following linking theorem:
\newpage
\begin{thm}[{\cite[Theorem 5.3]{R}}]\label{existence}
Let $E$ be a real Banach space with $E=V\oplus X$, where $V$ is
finite dimensional. Suppose ${\it\Phi}\in C^1(E, \mathbb{R})$,
satisfies $(PS)$-condition, and

$({\it\Phi}_1)$ there are constants $\rho, \alpha>0$ such that
${\it\Phi}\mid_{\partial B_{\rho}\cap X}\geq\alpha$, and

$({\it\Phi}_2)$ there is an $\;e\in \partial B_{\rho}\cap X$ and
$R>\rho$ such that if \;$Q\equiv (\overline{B}_R\cap V)\oplus
\{re\mid 0<r<R\}$, then ${\it\Phi}\mid_{\partial Q}\leq 0$.

\noindent where $B_r$ is an open ball in $E$ of radius $r$ centered
at $0$.

Then ${\it\Phi}$ possess a critical value $c\geq\alpha$ which can be
characterized as
\[
c\equiv \inf_{h\in\Gamma}\max_{u\in Q}{\it\Phi}(h(u)),
\]
where
\[
\Gamma=\{h\in C(\overline{Q}, E)\mid h={\rm id} \;on\; \partial Q\}.
\]
\end{thm}

The other one is the $\mathbb{Z}_2$-symmetric Mountain Pass Theorem:

\begin{thm}[{\cite[Corollary 7.22]{G}}]\label{multiplicity}
Let ${\it\Phi}$ be an even $C^1$-functional satisfying $(PS)$ on
$X=Y\oplus Z$ where $\dim (Y)=k<\infty$. Assume ${\it\Phi}(0)=0$ as
well as the following conditions:

$(1)$ There is $\rho>0$ and $\alpha\geq0$ such that $\inf
{\it\Phi}(S_{\rho}(Z))\geq\alpha$.

$(2)$ There exists $R>\rho$ and a subspace $F$ of $X$ containing $Y$
such that $\dim(F)=n>k$ and $\sup{\it\Phi}(S_R(F))\leq0$.

There exists then critical values $c_i$ $(1\le i\leq n-k)$ for
${\it\Phi}$ such that

$(a)$ $0\leq\alpha\leq c_1\leq\cdots\leq c_{n-k}$.

$(b)$ ${\it\Phi}$ has at least $n-k$ distinct pairs of non-trivial
critical points.
\end{thm}

In order to prove our main result by virtue of the above theorems,
we need to investigate the $(PS)$-condition and the linking
structure with respect to the functional. We will divide it into two
parts  and follow partially the ideas of the paper \cite{DJ} to give
the proofs of some lemmas in the two parts as follows.

\newpage
\noindent {\bf Part I. The $(PS)$-condition}

we will discuss the $(PS)$-condition in this part.

\begin{lem}\label{bounded}
Suppose that $\rm(W1), (W2)$ and $\rm(W4)$ are satisfied, %
then any $(PS)$-sequence is bounded.
\end{lem}

\pf Let $(u_n)\subset E$ be a $(PS)$-sequence, i.e., there exists a
constant $C_0>0$ such that
\begin{equation}\label{PSsequence}
  |{\it\Phi}(u_n)|\leq C_0\quad{\rm and}\quad
  {\it\Phi'}(u_n)\rightarrow 0.
\end{equation}

Arguing indirectly we assume that, up to a subsequence,
$\|u_n\|\rightarrow\infty$ and set $v_n=u_n/\|u_n\|$. Then
$\|v_n\|=1$. By Lemma \ref{cptem}, passing to a subsequence if
necessary, $v_n\rightharpoonup v$ in $E$ and $v_n\rightarrow v$ in
$L^p$ for all $2\leq p\leq\infty$. Then $v_n$ is bounded in
$L^\infty$. Since, by (W1) and (W2), $|w_u(t, u)|\leq C_w|u|$ for
some $C_w>0$, $w_u(t, u)=o(|u|)$ as $|u|\rightarrow\infty$, $\forall
t \in \mathbb{R}$ and $|u_n(t)|\rightarrow\infty$ if $v(t)\neq0$,
then it follows, by Lebesgue's Dominated Convergence Theorem, that

\begin{align*}
  \int\limits_{\mathbb{R}}\frac{\<W_u(t,u_n(t)), \varphi(t)\>}{\|u_n\|}dt%
  &=\int\limits_{\mathbb{R}}\<M(t)v_n(t), \varphi(t)\>dt%
  +\int\limits_{u_n(t)\neq
  0}\frac{\<w_u(t,u_n(t)), \varphi(t)\>|v_n(t)|}{|u_n|}dt\\
  &\rightarrow \int\limits_{\mathbb{R}}\<M(t)v(t), \varphi(t)\>dt\quad{\rm
  as}\; n\to\infty
\end{align*}
for all $\varphi\in C_0^{\infty}(\mathbb{R}, \mathbb{R}^{N})$.\\
By \eqref{derivative}, we have
\begin{equation}
\frac{{\it\Phi}'(u_n)\varphi}{\|u_n\|}%
=(v_n^+,\varphi)-(v_n^-,\varphi)%
-\int\limits_{\mathbb{R}}\frac{\<W_u(t,u_n(t)),\varphi(t)\>}{\|u_n\|}dt\notag
\end{equation}
for all $\varphi\in C_0^{\infty}(\mathbb{R}, \mathbb{R}^{N})$. From
this we deduce, using \eqref{PSsequence}, that
\[
  (-d^2/dt^2+L(t))v(t)=M(t)v(t),
\]
i.e.,
\begin{equation}\label{eigenvalue}
  (A-M)v=0.
\end{equation}
\par We claim that $v\neq0.$ Arguing by contradiction we assume that
$v=0.$ Choose $b>0$ in Lemma \ref{normcomp} such that
$\frac{C_W}{b}<1,$ where $C_W$ is the constant in (W1). Since
$E^{b-}\subset$ E in Lemma \ref{normcomp} is of finite-dimension,
then the compactness of the orthogonal projection
$P^{b-}:E\rightarrow E^{b-}\subset E$ implies $v_n^{b-}\rightarrow
v^{b-}=0$ in $E$.

It follows from \eqref{derivative} that
\[
  \frac{{\it\Phi'}(u_n)\left((u_n^{b+})^+-(u_n^{b+})^-\right)}{\|u_n\|^2}=\left\|v_n^{b+}\right\|^2-
  \int\limits_{u_n(t)\neq0}\frac{\<W_u(t,u_n),(v_n^{b+})^+-(v_n^{b+})^-\>}{|u_n|}|v_n|dt,
\]
then
\begin{align*}
 \left\|v_n^{b+}\right\|^2%
 &=\int\limits_{u_n(t)\neq0}\frac{\<W_u(t,u_n),(v_n^{b+})^+-(v_n^{b+})^-\>}{|u_n|}|v_n|dt%
 +\frac{{\it\Phi'}(u_n)\left((u_n^{b+})^+-(u_n^{b+})^-\right)}{\|u_n\|^2}\\[5pt]
 &\leq C_W\int\limits_{\mathbb{R}}\left|(v_n^{b+})^+-(v_n^{b+})^-\right||v_n|dt%
 +\frac{\|{\it\Phi'}(u_n)\|}{\|u_n\|}\\[7pt]%
 &\leq\frac{C_W}{2}\left(\int\limits_\mathbb{R}\left|(v_n^{b+})^++(v_n^{b+})^-\right|^2dt%
 +\int\limits_\mathbb{R}\left|(v_n^{b+})^+-(v_n^{b+})^-\right|^2dt\right)\\[7pt]%
 &\quad\quad+\frac{C_W}{2}\int\limits_\mathbb{R}\left|v_n^{b-}\right|^2dt%
 +\frac{\|{\it\Phi'}(u_n)\|}{\|u_n\|}\\[7pt]
 &= C_W \left|v_n^{b+}\right|_2^2+\frac{C_W}{2}\left|v_n^{b-}\right|_2^2%
 +\frac{\|{\it\Phi'}(u_n)\|}{\|u_n\|}\\[7pt]%
 &\leq\frac{C_W}{b}\left\|v_n^{b+}\right\|^2+\frac{C_W}{2}\left|v_n^{b-}\right|_2^2%
 +\frac{\|{\it\Phi'}(u_n)\|}{\|u_n\|},
\end{align*}
where $|\cdot|_2$ is the norm on $L^2$ and $(\cdot)^+$, $(\cdot)^-$
are the respective components with respect to the orthogonal
decomposition in Remark \ref{orthdecom}. The last inequality follows
by Lemma \ref{normcomp}. Note that $v_n^{b-}\rightarrow 0$ in $L^2$
since $v_n^{b-}\rightarrow$ in $E$. Thus
$\frac{C_W}{b}<1$ and \eqref{PSsequence} imply $\left\|v_n^{b+}\right\|^2\rightarrow 0.$ %
Then $1=\|v_n\|^2=\|v_n^{b-}\|^2+\|v_n^{b+}\|^2\rightarrow 0,$ a
contradiction.\par Therefore, $v\neq 0.$ Then \eqref{eigenvalue}
implies that $0$ is an eigenvalue of $A-M$ which is in contradiction
to (W4). $\hfill\Box$

\begin{lem}\label{PScodition}
Suppose that $\rm(W1),(W2)$ and $\rm(W4)$ are satisfied. Then
$\it\Phi$ satisfies the $(PS)$-condition.
\end{lem}
\pf Let $(u_n)\subset E$ be an arbitrary $(PS)$-sequence. By Lemma
\ref{bounded}, it is bounded, hence, we may assume without loss of
generality that $u_n\rightharpoonup u$ in $E$ and hence
$u_n^+\rightharpoonup u^+ $ and $u_n^-\to u^-$ due to
$\dim(E^-)<\infty$. By Lemma \ref{cptem},  $u_n\rightarrow u$ and
$u_n^+\to u^+$ in $L^2$. Observe that
\begin{align}\label{Cauchy}
\|u_n^+-u_m^+\|^2=&\left({\it\Phi}'(u_n)-{\it\Phi}'(u_m)\right)(u_n^+-u_m^+)\notag\\
&+\int\limits_{\mathbb{R}}\<W_u(t,u_n(t))-W_u(t,u_m(t)),
u_n^+-u_m^+\>dt,\;\forall n,m\in\mathbb{N}.
\end{align}
By (W1) and H\"{o}lder inequality
\begin{align*}
&\left|\int\limits_{\mathbb{R}}\<W_u(t,u_n(t))-W_u(t,u_m(t)),u_n^+-u_m^+\>dt\right|\\[10pt]
&\leq C_W\int\limits_{\mathbb{R}}(|u_n|+|u_m|)|u_n^+-u_m^+|dt\\[7pt]
&\leq C_W(|u_n|_2+|u_m|_2)|u_n^+-u_m^+|_2\to 0\quad {\rm as}\;n,m\to
\infty
\end{align*}
since $u_n\rightarrow u$ and $u_n^+\rightarrow u^+$ in $L^2$.\\
Note that
\[
\left({\it\Phi}'(u_n)-{\it\Phi}'(u_m)\right)(u_n^+-u_m^+)\to 0 \quad
{\rm as}\;n,m\to \infty
\]
since ${\it\Phi'}(u_n)\to 0$ and $(u_n)$ is bounded in $E$. Then
\eqref{Cauchy} implies that $(u_n^+)$ is a Cauchy sequence in $E$.
Hence $u_n^+\to u^+$ in $E$. Recall that $\dim(E^-\oplus
E^0)<\infty$, then $u_n^-+u_n^0\to u^-+u^0$ in $E$. This yields
$u_n\to u$ in $E$ and the proof is completed. $\hfill\Box$

\noindent{\bf Part II. Linking structure.}

First we have the following lemma.
\begin{lem}\label{linksphere}
Let $\rm(W1)$ be satisfied. Then there exists $\rho>0$ such that
\[
\alpha:=\inf{\it\Phi}(\partial B_{\rho}\cap E^+)>0.
\]
\end{lem}

\pf By Lemma \ref{cptem}, we have
\begin{equation}\label{normcontrol1}
|u|_{\infty}\rightarrow 0 \;{\rm as} \;\|u\|\rightarrow 0,
\end{equation}
where $|\cdot|_{\infty}$ is the norm on $L^{\infty}.$ From (W1), we
obtain that $W(t,u)=o(|u|^2)$ as $|u|\rightarrow 0$ uniformly in
$t$. Combining this with \eqref{normcontrol1}, for any
$\varepsilon>0,$ there is a $\delta>0$ such that
\[
{\it\Psi}(u)\leq\varepsilon |u|_{2}^{2}\leq\varepsilon
\beta_2^2\|u\|^{2},\; \forall \;\|u\|\leq \delta,
\]
where $\beta_2$ is the constant in \eqref{LpE}. Taking
$\varepsilon=1/(4\beta_2^2)$ and $0<\rho<\delta$, then
$\alpha:=\inf{\it\Phi}(\partial B_{\rho}\cap E^+)\geq \rho^2/4>0$ by
the form of $\it\Phi$ in \eqref{functional}. $\hfill\Box$

Due to (W3) and the spectral result of $A$ in the previous section,
we can arrange all the eigenvalues (counted with multiplicity) of
$A$ in $(0, m_0)$ by $0<\lambda_1\leq \lambda_2\leq \cdots\leq
\lambda_{\ell}<m_0$ and let $e_j$ denote the corresponding
eigenfunctions: $Ae_j=\lambda_je_j$ for $j=1, \ldots, \ell$. Set
$E_{\ell}^+:=\mathrm{span}\{e_1, \ldots, e_{\ell}\}$. According to
the definition of the norm on $E$, we have
\begin{equation}\label{normcontrol2}
\lambda_1|v|_2^2\leq \|v\|^2\leq \lambda_{\ell}|v|_2^2\quad\mbox{for
all }v\in E_{\ell}^+.
\end{equation}
Set $\tilde{E}=E^-\oplus E^0\oplus E_{\ell}^+$.

\begin{lem}\label{linkQ}
Let $\rm(W1),(W2)$ and $\rm(W3)$ be satisfied and $\rho>0$ be given
by Lemma $\ref{linksphere}$. Then there exists $R_{\tilde{E}}>\rho$
such that ${\it\Phi}(u)<0$ for all $u\in \tilde{E}$ with $\|u\|\geq
R_{\tilde{E}}$.
\end{lem}

\pf It suffice to show that ${\it\Phi}(u)\rightarrow -\infty$ as
$u\in \tilde{E}$, $\|u\|\rightarrow \infty.$ Arguing indirectly we
assume that there exist some $c>0$ and a sequence $(u_j)\subset
\tilde{E}$ with $\|u_j\|\rightarrow \infty$ such that
${\it\Phi}(u_n)\geq -c$ for all $n$. Then, setting
$v_n=u_n/\|u_n\|,$ we have $\|v_n\|=1$, and we may assume without
loss of generality $v_n\rightarrow v$, $v_n^-\rightarrow v^-,
v_n^0\rightarrow v^0, v_n^+\rightarrow v^+\in E_{\ell}^+$ since
$\dim (\tilde{E})<\infty$.

From \eqref{functional}, we have
\begin{equation}\label{contrineq1}
-\frac{c}{\|u_n\|^2}\leq \frac{{\it\Phi}(u_n)}{\|u_n\|^2}=
\frac{1}{2}\|v_n^+\|^2-\frac{1}{2}\|v_n^-
\|^2-\int\limits_{\mathbb{R}}\frac{W(t,u_n)}{\|u_n\|^2}dt.
\end{equation}
We claim that $v^+\neq0.$ Indeed, if not it follows from
\eqref{contrineq1} and (W1) that $\|v_n^-\|\rightarrow 0$ and thus
$v_n\rightarrow v=v^0$.  Also
$\int\limits_{\mathbb{R}}\frac{W(t,u_n)}{\|u_n\|^2}dt\rightarrow 0$
and, by Lemma \ref{cptem}, $v_n\rightarrow v$ in $L^2$.

Note that by (W1) and (W2), $W(t, u)=\frac{1}{2}M(t)u\cdot u +w(t,
u)$ and $|w(t, u)|\leq C_w|u|^2$ for some $C_w>0$, $w(t,
u)/|u|^2\rightarrow 0$ as $|u|\rightarrow \infty$, $\forall
t\in\mathbb{R}$. Since $|u_n(t)|\rightarrow \infty$ if $v(t)\neq 0$,
we obtain
\begin{align}\label{contrineq2}
\int\limits_{\mathbb{R}}\frac{|w(t,u_n)|}{\|u_n\|^2}dt
&=\int\limits_{u_n(t)\neq 0}\frac{|w(t,u_n)|}{|u_n|^2}|v_n|^2dt\notag\\
&\leq2\int\limits_{u_n(t)\neq 0}\frac{|w(t, u_n)|}{|u_n|^2}|v_n-v|^2dt+%
2\int\limits_{u_n(t)\neq 0}\frac{|w(t, u_n)|}{|u_n|^2}|v|^2dt\notag\\[5pt]
&\leq 2C_w\int\limits_{\mathbb{R}}|v_n-v|^2dt+
2\int\limits_{u_n(t)\neq 0}\frac{|w(t, u_n)|}{|u_n|^2}|v|^2dt\notag\\
&=o(1).
\end{align}
The last equality holds by $v_n\rightarrow v$ in $L^2$ and
Lebesgue's Dominated Convergence Theorem. Also, by (W3),
\begin{equation}\label{contrineq3}
\frac{1}{2}\int\limits_{\mathbb{R}}\frac{\<M(t)u_n,
u_n\>}{\|u_n\|^2}dt= \frac{1}{2}\int\limits_{u_n(t)\neq
0}\frac{\<M(t)u_n,
u_n\>}{|u_n|^2}|v_n|^2dt\geq\frac{m_0}{2}|v_n|_2^2.
\end{equation}

\noindent From \eqref{contrineq2}, \eqref{contrineq3} and since
$\int\limits_{\mathbb{R}}\frac{W(t,u_n)}{\|u_n\|^2}dt\rightarrow 0$
it follows that $|v_n|_2\rightarrow 0$. Due to $\dim
(\tilde{E})<\infty$, $1=\|v_n\|\rightarrow 0$ and this contradiction
implies that $v^+\neq 0$. Note that (W3), \eqref{normcontrol2} and
Remark \ref{orthdecom} imply that
\begin{align*}
  \|v^+\|^2-\|v^-\|^2-\int\limits_{\mathbb{R}}\<M(t)v, v\>dt
   &\leq\|v^+\|^2-\|v^-\|^2-m_0|v|_2^2\\
   &\leq-\left((m_0-\lambda_{\ell})|v^+|_2^2+\|v^-\|^2+m_0|v^-+v^0|_2^2\right)<0
\end{align*}

\noindent Then there is $T>0$ such that
\begin{equation}\label{contrineq4}
\|v^+\|^2-\|v^-\|^2-\int\limits_{-T}^T\<M(t)v, v\>dt<0.
\end{equation}
By \eqref{contrineq2}, we get
\[
  \lim_{n\rightarrow \infty}\int\limits_{-T}^T\frac{w(t,u_n)}{\|u_n\|^2}dt
  \to 0.
\]

\noindent Thus \eqref{contrineq1} and \eqref{contrineq4} imply that
\begin{align*}
  0&\leq\lim_{n\rightarrow \infty}\left(\frac{1}{2}\|v_n^+\|^2-\frac{1}{2}\|v_n^-\|^2
  -\int\limits_{-T}^T\frac{W(t,u_n)}{\|u_n\|^2}dt\right)\\
  &=\frac{1}{2}\left(\|v^+\|^2-\|v^-\|^2-\int\limits_{-T}^T\<M(t)v, v\>dt\right)<0,
\end{align*}
a contradiction. $\hfill\Box$

As an immediate result of Lemma \ref{linkQ}, we have
\begin{lem}\label{sublinkQ}
Let $\rm(W1),(W2)$ be satisfied and $\rho>0$ be given by Lemma
$\ref{linksphere}$. Then, letting $e\in E_{\ell}^+$ with
$\|e\|=1,$ there is $R>\rho$ %
such that  $\sup{\it\Phi(\partial Q)}\leq0$ where
$Q:=\{u=u_1+re:u_1\in E^-\oplus E^0, \;\|u_1\|\leq R,\; 0<r<R\}.$
\end{lem}
\pf Set $R=R_{\tilde{E}}$, where $R_{\tilde{E}}$ is the constant in
Lemma \ref{linkQ}. Then, by Lemma \ref{linkQ},
\begin{equation}\label{Phi-e}
{\it\Phi}(u)<0, \; \forall u\in E^-\oplus
E^0\oplus\mathrm{span}\{e\}\subset \tilde{E}, \;\|u\|\geq R.
\end{equation}
Observe that
\[
\partial Q=Q_1\cup Q_2\cup Q_3,
\]
where
\begin{align*}
Q_1&:=\{u\in E^-\oplus E^0:\|u\|\leq R\},\\
Q_2&:=\{u=u_1+Re:u_1\in E^-\oplus E^0,\;\|u_1\|\leq R\},\\
Q_3&:=\{u=u_1+re:u_1\in E^-\oplus E^0,\;\|u_1\|=R,\;0\leq r\leq R\}.
\end{align*}
Due to \eqref{Phi-e}, it holds that
\[
{\it\Phi}(u)\leq 0,\;\forall u\in Q_2\cup Q_3.
\]
Also, in view of (W1) and the form of ${\it\Phi}$ in
\eqref{functional}, ${\it\Phi}(u)\leq 0,\;\forall u\in Q_1$. Then
the proof is completed. $\hfill\Box$

After all the above preparations, we now come to the proof of our
main result.

\noindent {\bf Proof of Theorem \ref{mainresult}.} Step 1.
Existence. With $V=E^-\oplus E^0$ and $X=E^+$ in Theorem
\ref{existence}, the conditions $({\it\Phi}_1)$ and $({\it\Phi}_2)$
there hold by Lemmas \ref{linksphere} and \ref{sublinkQ}
respectively. ${\it\Phi}$ satisfies  the $(PS)$-condition by Lemma
\ref{PScodition}. Hence, ${\it\Phi}$ has at least one critical point
$u$ with ${\it\Phi}(u)\geq \alpha>0$ by Theorem \ref{existence}.
Since ${\it\Phi}(0)=0$, $u$ is a nontrivial critical point of
${\it\Phi}$. Then \eqref{HS} has at least one nontrivial homoclinic
solution $u$ by Proposition \ref{propofPhi}.

Step 2. Multiplicity. Let $X=E$, $Y=E^-\oplus E^0$ and $Z=E^+$ in
Theorem \ref{multiplicity}. Since $W(t, u)$ is even in u, then
${\it\Phi}$ is even and ${\it\Phi}(0)=0$ by the form of ${\it\Phi}$
in \eqref{functional}. Lemma \ref{linksphere} shows that $(1)$ in
Theorem \ref{multiplicity} holds. With $F=\tilde{E}$ in Theorem
\ref{multiplicity}, then Lemma \ref{linkQ} implies that $(2)$ in
Theorem \ref{multiplicity} also holds. Note that
$\dim(F)-\dim(Y)=\dim(E_\ell^+)=\ell$. Therefore, ${\it\Phi}$ has at
least $\ell$ pairs of nontrivial critical points by Theorem
\ref{multiplicity} and then \eqref{HS} has at least $\ell$ pairs of
nontrivial homoclinic solutions  by Proposition \ref{propofPhi}.
$\hfill\Box$


\newpage
\begin{thebibliography}{99}\small
\item F. Antonacci, Periodic and homoclinic solutions to a class of
Hamiltonian systems with indefinite potential in sign, Boll. Un.
Mat. Ital. B (7) 10 (1996) 303--324.
\item P.C. Carriao, O.H. Miyagaki, Existence of homoclinic
solutions for a class of time-dependent Hamiltonian systems, J.
Math. Anal. Appl. 230 (1999) 157--172.
\item C.N. Chen, S.Y. Tzeng, Existence and multiplicity results for homoclinic orbits of
Hamiltonian systems, Electron. J. Differential Equations 1997 (1997)
1--19.
\item V. Coti Zelati, P.H. Rabinowitz, Homoclinic orbits for
second order Hamiltonian systems possessing superquadratic
potentials, J. Amer. Math. Soc. 4 (1991) 693--727.
\bibitem{D} Y.H. Ding, Existence and multiplicity results for homoclinic
solutions to a class of Hamiltonian systems, Nonlinear Anal. 25
(1995) 1095--1113.
\item Y.H. Ding, M. Girardi, Periodic and homoclinic solutions
to a class of Hamiltonian systems with the potentials changing sign,
Dynam. Systems Appl. 2 (1993) 131--145.
\bibitem{DJ} Y.H. Ding, L. Jeanjean, Homoclinic orbits for a nonperiodic Hamiltonian system, J.
Differential Equations 237 (2007) 473--490.
\item G.H. Fei, The existence of homoclinic orbits for
Hamiltonian sytems with the potential changing sign, Chinese Ann.
Math. Ser. B 17 (1996) 403--410.
\item P.L. Felmer, E.A. De B.e. Silva, Homoclinic and periodic orbits for
Hamiltonian systems, Ann. Sc. Norm. Super. Pisa Cl. Sci. (4) 26 (2)
(1998) 285--301.
\bibitem{G} N. Ghoussoub, Duality and Perturbation Methods in
Critical Point Theory, Cambridge University Press, Cambridge, 1993.
\bibitem{K} T. Kato, Perturbation Theory for Linear Operators,
Springer-Verlag, New York, 1980.
\item P. Korman, A.C. Lazer, Homoclinic orbits for a class of
symmetric Hamiltonian systems, Electron. J. Differential Equations
1994 (1994) 1--10.
\bibitem{R} P.H. Rabinowitz, Minimax Methods in Critical Point Theory with
Applications to Differential Equations, in: CBMS Regional Conf. Ser.
in Math., vol. 65, American Mathematical Society, Providence, RI,
1986.
\item Z.Q. Ou, C.L. Tang, Existence of homoclinic solutions for the second order Hamiltonian
systems, J. Math. Anal. Appl. 291 (1) (2004) 203--213.
\item P.H. Rabinowitz, Homoclinic orbits for a class of Hamiltonian
systems, Proc. Roy. Soc. Edinburgh Sect. A 114 (1990) 33--38.
\item P.H. Rabinowitz, K. Tanaka, Some results on connecting
orbits for a class of Hamiltonian systems, Math. Z. 206 (1990)
473--499.
\item A. Salvatore, Homoclinic orbits for a special class of nonautonomous Hamiltonian
systems, in: Proceedings of the Second World Congress of Nonlinear
Analysts, Part 8, (Athens, 1996), Nonlinear Anal. 30 (8) (1997)
4849--4857.
\item S.P. Wu, J.Q. Liu, Homoclinic orbits for second order Hamiltonian
system with quadratic growth, Appl. Math. J. Chinese Univ. Ser. B 10
(1995) 399--410.
\item S.P. Wu, H.T. Yang, A note on homoclinic orbits for second order Hamiltonian
system, Appl. Math. J. Chinese Univ. Ser. B 13 (1998)  251--262.

\end{thebibliography}
\end{document}